\documentclass[a4paper,10pt]{article}

\def \1{{\bf 1}}
\def \N {{\mathbf {N}}}
\def \K {{\cal K}}
\def \eps {{\varepsilon}}
\title{On Thouvenot's ergodic proof of Roth's theorem }
\author{V.V. Ryzhikov}

\begin{document}
\Large
\maketitle

%\begin{abstract}{\large    a bit later}
%\end{abstract}
\section{Introduction}
Roth's theorem says that a subset of $\N$ of a  positive density contains an arithmetic progression of length 3.
H. Furstenberg has  proved that this theorem  is equivalent to the following assertion: 
\it for any invertible measure-preserving transformation $T$ of a probability space $(X,\mu)$ and
any set $A$ of a positive measure it holds 
$$\liminf_N\frac{1}{N}\sum_{i=1}^{N}\mu(A\cap T^iA\cap T^{2i}A)>0.
$$
\rm Furstenberg  also gave an ergodic proof of this fact (see \cite{F}).  In 2002 J.-P. Thouvenot communicated me an interesting  modification
of this  proof using an observation from \cite{R} (see also a joining proof of Marcus' theorem on multiple mixing for horocycle flows \cite{T}) .  Sometimes I included his short proof  in my talks replacing a  joining by an operator.
 Now I present this topic here (section 2) 
adding  old remarks-proofs connected with  Furstenberg's theorems on multiple progression average mixing for weakly mixing transformations (section 3).

\section{Thouvenot's  proof of Furstenberg's version of Roth's theorem }

Let $f,g,h\in L_\infty(X,\mu)$.  From any sequence $N_{k'}$ we choose a subsequence $N_k$ such that  for an operator $J:L_2\to L_2\otimes L_2$ the equality
$$\langle Jf, g\otimes h\rangle_{L_2\otimes L_2} =\lim_{k}\frac{1}{N_k}\sum_{i=1}^{N_k} \int f T^igT^{2i}h \ d\mu$$
holds for any
$f,g,h\in L_\infty(X,\mu)$. 
The  definition of $J$ is correct, this follows from the ergodicity of $T$:
$$\lim_{k}\frac{1}{N_k}\sum_{i=1}^{N_k} \int  T^igT^{2i}h \ d\mu =\lim_{k}\frac{1}{N_k}\sum_{i=1}^{N_k} \int  gT^{i}h \ d\mu =
\int f d\mu \int gd\mu.$$

We see that $(T\otimes T^2)J=J$.   So $(T\otimes T^2)Jf=Jf$.  But a $(T\otimes T^2)$-invariant
function belongs to $L_2(\K\otimes \K, \mu\otimes\mu)$, where $\K$ is a compact factor algebra ( = Kronecker  algebra generated by all proper functions of $T$).  Indeed, we must only  to remark that  a restriction of $T$ (and $T^2$ as well)  onto  $L_2(\K,\mu)^\perp$,  say $T'$, has the  property   $T'^{i}\to_w 0$ ($T'$ has continuous spectrum),  hence, $(T\otimes T^2)F=F$ implies $F\in    L_2(\K,\mu)    \otimes  L_2(\K,\mu) $.                                                             

 Denoting   $P$ for the orthogonal  projection $L_2(\mu)\to L_2(\K,\mu)$
we obtain 
$$\langle Jf, g\otimes h\rangle =\langle Jf, Pg\otimes Ph\rangle =
\lim_{k}\frac{1}{N_k}\sum_{i=1}^{N_k} \int f \ T^iPg\ T^{2i}Ph\  d\mu =$$
$$=\lim_{k}\frac{1}{N_k}\sum_{i=1}^{N_k} \int Pf \ T^iPg \ T^{2i}Ph \ d\mu.$$
\rm
Let $f=g=h=\chi_A$, $\mu(A)>0$.
A closure of $\{T^iPf \}$ is a compact set, for any $\eps>0$ there is $L$ such that for any $n$ and for at least one of
 $i= n+1, n+2\dots, n+L$ we get 
$$\|T^iPf-Pf\|_{L_2}<\eps.$$
Thus, for a sufficiently small $\eps'>0$ we have
$$\liminf_N\frac{1}{N}\sum_{i=1}^{N}\mu(A\cap T^iA\cap T^{2i}A)\geq \frac{1}{L}\left(\int (P\chi_A)^3\ d\mu \ -\eps'\right)>0.
$$
\section{ Remarks to  Furstenberg's theorems on weakly mixing transformations}
    Furstenberg   \cite{F}  proved the following theorem:  If $T$ is weakly mixing, then 
$$\frac{1}{N}\sum_{i=1}^{N} \prod_{p=1}^{m}T^{pi}f_p\ \to_{L_2}\ \prod_{p=1}^{m}\int f_p \ d\mu \ \ (N\to\infty)
\eqno (2, m)$$
holds for any collection of $f_i\in L_\infty$.
Let $f,g,h\in L_\infty(X,\mu)$ and $T$ be weakly mixing,  let us show  
$$\frac{1}{N}\sum_{i=1}^{N} \int f T^{i}gT^{2i}h \to \ \int f d\mu \int gd\mu\int hd\mu, \eqno (1,2)$$
$$\frac{1}{N}\sum_{i=1}^{N}  T^if T^{2i}gT^{3i}h \ \to_{L_2} \ \int f d\mu \int gd\mu\int hd\mu.\eqno(2,3)$$

\bf Proof of (1,2). \rm We define a joining 
$$\nu(f\otimes g\otimes h)=\lim_{k} \frac{1}{N_k}\sum_{i=1}^{N_k} \int f T^{i}gT^{2i}h.$$
We have $(I\otimes T\otimes T^2)\nu=\nu$, but $T\otimes T^2$ is ergodic.
$$\nu(f\otimes g\otimes h)=\nu\left(f\otimes \left( \frac{1}{N}\sum_{i=1}^{N}  T^{i}g\otimes T^{2i}h\right)\right)=$$
$$\nu(f\otimes \1\otimes \1)\int gd\mu\int hd\mu =\int f d\mu \int gd\mu\int hd\mu.$$ 
 (1) is proved. Here we can use also that $Id$ and an ergodic transformation $S=T\otimes T^2$ are disjoint, so our joining has to be a direct product of its projections,  see \cite{R},
\cite{T}. 

\bf Proof of (2,3). \rm  We define a joining $\eta$ setting 
$$\eta( f\otimes g\otimes h\otimes f'\otimes g'\otimes h')=\lim_k\frac{1}{N_k^2}\int\sum_{i=1}^{N_k} T^if \  T^{2i}g\ T^{3i}h\ \
\sum_{j=1}^{N_k} \ T^j f' \ T^{2j}g'\ T^{3j}h'\  d\mu.$$
From the above definition it follows an invariance
$$\eta=(I\otimes I \otimes I\otimes T\otimes T^2\otimes T^3)\eta,$$
but $T\otimes T^2\otimes T^3$ is ergodic. Again our joining will be a product: $\eta=\mu^3\otimes\mu^3$.
Here we have made use  of (1,2):  the  projections of $\eta$ are equal to $\mu^3$, indeed 
$$\eta(f\otimes g\otimes h\otimes \1\otimes\1\otimes\1)=\lim_N\frac{1}{N}\sum_{i=1}^N\int T^if\ T^{2i}g\ T^{3i}h \ d\mu=$$
$$= \lim_N\frac{1}{N}\sum_{i=1}^N\int f\ T^{i}g\ T^{2i}h \ d\mu=  \ \int f d\mu \int gd\mu\int h\ d\mu.$$
Let   $\int f d\mu=0$,  then 
$$
\lim_{k}\|\frac{1}{N_k}\sum_{i=1}^{N_k}  T^if T^{2i}gT^{3i}h\|^2_{L_2}=
\eta( f\otimes g\otimes h\otimes  f\otimes g\otimes h)=0 .\eqno (2',3)
$$
To prove  (2,3) we have   to say only that for any sequence $N_{k'}$ one can choose a subsequence $N_{k}$ for
which (2',3) holds. 

Now let's remark that (2,3) implies 
$$\frac{1}{N}\sum_{i=1}^{N} \int f_0 T^{i}f_1\ T^{2i}f_2 \ T^{3i}f_3 \ d\mu\ \to \prod_{p=0}^3  \int f_p \ d\mu \ \ (N\to\infty).  \eqno (1,3)$$
From (1-3) we deduce  as above
$$\frac{1}{N}\sum_{i=1}^{N}  T^if_1 \ T^{2i}f_2\ T^{3i}f_3 \ T^{4i}f_4\ \to_{L_2} \   \prod_{p=1}^4  \int f_p \ d\mu, \eqno (2,4)$$
and so on:  (2,m) implies (1, m), from (1,m) we get 
 $$\frac{1}{N}\sum_{i=1}^{N} \prod_{p=1}^{m+1}T^{pi}f_p\ \to_{L_2}\ \prod_{p=1}^{m+1}\int f_p \ d\mu
\eqno (2, m+1)$$
as $N\to\infty$.

vryzh@mail.ru
\end{document}